\documentclass[a4paper,12pt]{amsart}

\usepackage{graphicx}
\usepackage{amsmath}

\textwidth=\paperwidth \advance\textwidth by-2in
\makeatletter
\def\temp{12}
\ifx\@mainsize\temp \headheight=6.5pt \fi
\makeatother

\def\Z{\mathbf{Z}}
\def\Q{\mathbf{Q}}
\def\R{\mathbf{R}}
\def\C{\mathbf{C}}
\def\L{\mathbf{L}}
\def\S{\mathbf{S}}
\def\sign{\operatorname{sign}}
\def\sgn{\operatorname{sgn}}
\def\rank{\operatorname{rank}}
\def\lk{\operatorname{lk}}
\def\Ker{\operatorname{Ker}}

\def\inte{\operatorname{int}}
\def\bd{\partial}
\def\spmatrix#1{\left[
  \begin{smallmatrix} #1 \end{smallmatrix}
  \right]\ignorespaces}
\theoremstyle{plain}
\newtheorem{thm}{Theorem}[section]

\newtheorem{lem}{Lemma}[section]
\theoremstyle{definition}
\newtheorem*{exmp}{Example}
\theoremstyle{remark}

\begin{document}

\title{Signatures of Links in Rational Homology Spheres}
\date{August 28, 2001 (First Edition: November 23, 1999)}
\author{Jae Choon Cha}
\email{jccha\char`\@knot.kaist.ac.kr}
\address{Department of Mathematics\\
Korea Advanced Institute of Science and Technology\\
Taejon, 305--701\\
Korea}
\author{Ki Hyoung Ko}
\email{knot\char`\@knot.kaist.ac.kr}
\address{Department of Mathematics\\
Korea Advanced Institute of Science and Technology\\
Taejon, 305--701\\
Korea}

\begin{abstract}
\sloppy A theory of signatures for odd-dimensional links in rational
homology spheres is studied via their generalized Seifert
surfaces. The jump functions of signatures are shown invariant under
appropriately generalized concordance and a special care is given to
accommodate 1-dimensional links with mutual linking.  Furthermore our
concordance theory of links in rational homology spheres remains highly
nontrivial after factoring out the contribution from links in integral
homology spheres.
\end{abstract}

\maketitle

\section{Introduction and main results}

Signature invariants of knots and links have been studied in many
articles, and turned out to be effective in various problems in knot
theory, especially in studying concordance.  Signatures were first
defined as numerical invariants from Seifert matrices of classical
knots in the three space by Trotter~\cite{Tro}.  Murasugi formulated a
definition of signatures for links~\cite{Mu}.  Tristram introduced
link signature functions~\cite{Tr}.  Kauffman and Taylor found a
geometric approach to signatures of links using branched
covers~\cite{KT}.  Milnor defined a signature function from a duality
in the infinite cyclic cover of a knot complement~\cite{M1}. It is
known that all the approaches are
equivalent~\cite{Er,KT,Ma}. Signature invariants are known to be
concordance invariants of codimension two knots and links in odd
dimensional spheres. In the remarkable works of Levine~\cite{L1,L2},
signature invariants are used to classify concordance classes of
knots.  Levine showed that the knot concordance groups are isomorphic
to the cobordism groups of Seifert matrices in higher odd dimensions,
and signatures are complete invariants for the matrix cobordism groups
modulo torsion elements.

In this paper, we develop signature invariants of arbitrary links of
codimension two in odd dimensional rational homology spheres,
motivated from the work of Cochran and Orr~\cite{CO1,CO2}.  Our
signatures are invariant under the following equivalence of links in
rational homology spheres. An $n$-\emph{link} $L$ in a rational
homology $(n+2)$-sphere $\Sigma$ is a submanifold diffeomorphic to an
ordered union of disjoint standard $n$-spheres.  Such a link is
sometimes denoted by a pair $(\Sigma,L)$. Two $n$-links $(\Sigma_0,
L_0)$ and $(\Sigma_1, L_1)$ are \emph{rational homology concordant}
(or \emph{$\Q$-concordant}) if there is an $(n+3)$-manifold $\S$ whose
boundary is the disjoint union of $-\Sigma_0$ and $\Sigma_1$ such that
the triad $(\S, \Sigma_0, \Sigma_1)$ has the rational homology of
$(S^{n+2}\times [0,1], S^{n+2}\times 0, S^{n+2}\times 1)$, and there
is a proper submanifold $\L$ in $\S$ diffeomorphic to $L_0\times
[0,1]$ such that the boundary of the $i$-th component of $\L$ is the
union of the $i$-th components of $-L_0$ and $L_1$.  Concordant links
are clearly rational homology concordant as well, and so the
signatures studied in this paper are invariant under ordinary
concordance as well.

As done in~\cite{Tro,Mu,Tr,KT} for links in~$S^{2q+1}$, we extract
invariants of links from the signature function
$$
\sigma^q_A(\phi)=\sign((1-e^{i\phi})A+\epsilon(1-e^{-i\phi})\overline{A}^T)
$$
and its jump function
$$ \delta^q_A(\theta)=\begin{cases}
\lim\limits_{\phi\to\theta^+}\sigma^q_A(\phi)
-\lim\limits_{\phi\to\theta^-}\sigma^q_A(\phi), & \text{if
$\theta$ is not a multiple of $2\pi$} \\ \,\,0, & \text{otherwise}
\end{cases}
$$
where $A$ is a rational Seifert matrix, $\phi$ and $\theta$ are reals,
and $\epsilon=(-1)^{q+1}$.  We follow the convention that the
signature of an anti-hermitian matrix $P$ is the signature of the
hermitian matrix $\sqrt{-1}\cdot P$.
Section~\ref{sec:matrix-sign-poly} is devoted to an algebraic study of
$\delta^q_A$ and $\sigma^q_A$. Especially, we prove a
\emph{reparametrization formula} for $\delta^q_A$ and $\sigma^q_A$ of
a Seifert matrix of a union of parallel copies of a Seifert surface
(see Lemma~\ref{lem:alg-sign-formula},
Theorem~\ref{thm:sign-formula}).

A Seifert surface of a link $L$ is defined to be an oriented
submanifold bounded by $L$ in the ambient space. $L$ has a Seifert
surface if and only if there exists a homomorphism of the first
integral homology of the exterior of $L$ onto the infinite cyclic
group sending each meridian of $L$ to a fixed generator. In general,
such a homomorphism does not exist and links bound no Seifert
surfaces, unless the ambient space is an integral homology sphere.  In
the case of rational homology spheres, we show that a certain union of
parallel copies of a link $L$ bounds a Seifert surface even though $L$
itself does not.  More precisely, for an $m$-component link $L$ and
$\tau=(\tau_i)\in \Z^m$ where $\tau_i$ are relative prime, we call an
oriented submanifold $F$ a \emph{generalized Seifert surface of type
$\tau$ for $L$} if for some nonzero integer $c$, $F$ is bounded by the
union of $|c\tau_i|$ parallel copies (or negatively oriented copies if
$c\tau_i<0$) of the $i$-th component of $L$ for all $i$. If $L$ is a
1-link, the parallel copies are taken with respect to some framing $f$
on $L$ and the surface $F$ and the framing $f$ must induce the same
framing on the parallel copies.  We call $c$ the \emph{complexity
of~$F$}. A higher dimensional link admits a generalized Seifert
surface of any type $\tau$.  But a 1-link admits a generalized Seifert
surface of the types satisfying a certain condition and a framing is
uniquely determined by a given type.  Furthermore $\Q$-concordant
links admit generalized Seifert surfaces of the same type.  (See
Theorem~\ref{thm:seifert-surface-existence-highdim} and
Theorem~\ref{thm:seifert-surface-existence-dim1}.)

We define a Seifert matrix to be a matrix associated to the
rational-valued Seifert pairing on the middle dimensional rational
homology of a generalized Seifert surface, and derive invariants of a
link from the signature jump function of a Seifert matrix as follows.
Let $\sgn a = a/|a|$ for $a\ne 0$.  The \emph{signature jump function
of type $\tau$} for a $(2q-1)$-link $L$ is defined by
$$
\delta^\tau_L(\theta) = \sgn c \cdot \delta^q_A(\theta/c)
$$
where $A$ is a Seifert matrix of a generalized Seifert surface $F$ of
type $\tau$ for $L$ and $c$ is the complexity of~$F$.  Our main result
is that the signature jump function $\delta^\tau_L$ is a well-defined
invariant under $\Q$-concordance which is independent of the choice of
generalized Seifert surfaces.

\begin{thm}\label{thm:link-sign-inv}
If $L_0$ and $L_1$ are $\Q$-concordant, $\delta^\tau_{L_0}(\theta) =
\delta^\tau_{L_1}(\theta)$ for all~$\theta$.
\end{thm}

For $\tau=(1\cdots1)$, we denote $\delta^\tau_L$
by~$\delta^{\vphantom{\tau}}_L$.  If $L$ is a link in an integral
homology sphere, $L$ bounds a Seifert surface $F$ and
$\delta_{L}(\theta)$ can be computed from a Seifert matrix of~$F$.
Hence our signature invariant is a generalization of the well-known
signature invariant of links in the spheres.

For 1-links that are null-homologous modulo $d$ in a rational homology
3-sphere, Gilmer defined the $d$-signature that depends on the choice
of a specific element in the first homology of the ambient
space~\cite{G}. We remark that it is equal to the evaluation of the
signature function associated to a generalized Seifert matrix of type
$(1\cdots 1)$ at the $d$-th roots of unity, and so it can be derived
from our signature jump function.

A link is called a \emph{slice link} if its components bound properly
embedded disjoint disks in a rational homology ball bounded by the
ambient space.  For a link $L$ and an $n$-tuple
$\alpha=(s_1,\ldots,s_n)$ with $s_i=\pm 1$, let $i_\alpha L$ be the
union of $n$ parallel copies of $L$ where the $i$-th copy is oriented
according to the sign of~$s_i$, and let $n_\alpha$ be the sum
of~$s_i$.  The \emph{complexity $c^\tau(L)$ of type $\tau$} for $L$ is
defined to be the minimum of the absolute values of complexities of
generalized Seifert surfaces of type~$\tau$.

For $i=1,2$, let $(\Sigma_i,L_i)$ be an $n$-dimensional link with
$m_i$-components and $B_i$ be an $(n+2)$-ball in $\Sigma_i$ such that
for some $k \ge 1$, $B_i$ intersects $k$ components of $L_i$ at a
trivial $k$-component disk link in $B_i$.  By pasting $(\Sigma_1-\inte
B_1, L_1-\inte B_1 \cap L_1)$ and $(\Sigma_2-\inte B_2, L_2-\inte B_2
\cap L_2)$ together along an orientation reversing homeomorphism on
boundaries, a $(m_1+m_2-k)$-component link is obtained and is called a
\emph{connected sum of $L_1$ and $L_2$}. A connected sum depends on
the choice of $B_i$ if~$m\ge2$.  Some useful properties of our
signature invariant are listed in the following theorem.

\begin{thm}\label{thm:link-sign-prop}
\mbox{}
\begin{enumerate}
\item
If $L$ is slice, $\delta^\tau_L(\theta)=0$.
\item
If $n_\alpha\ne 0$, $\delta^\tau_{i_\alpha L}(\theta) = \sgn n_\alpha \cdot
\delta^\tau_L(n_\alpha\theta)$.
\item
$\delta^\tau_L(\theta)$ is a periodic function of period $2\pi\cdot
c^\tau(L)$.
\item
Let $L$ be a connected sum of $(2q-1)$-links $L_1$ and $L_2$.  If either
$q>1$ or $L_1$, $L_2$ are knots, then $\delta^\tau_L(\theta) =
\delta^\tau_{L_1}(\theta) + \delta^\tau_{L_2}(\theta)$.
\end{enumerate}
\end{thm}

In Section~\ref{sec:sign-of-link}, we prove the above results by
combining the techniques of signatures of links in $S^{2q+1}$ with the
reparametrization formula in Section~\ref{sec:matrix-sign-poly} that
plays a key role in resolving the difficulty that links do not bound
Seifert surfaces.

We give examples of links in rational homology spheres whose signature
jump functions have arbitrary rational multiples of $2\pi$ as periods,
and show that the theory of $\Q$-concordance of links in rational
homology spheres is not reduced to the concordance theory of links in
spheres.  In fact we construct a family of infinitely many
$\Q$-concordance classes of knots that are independent modulo knots in
the spheres in the following sense.  Let $C_n$ and $C^\Q_n$ be the set
of concordance classes of $n$-knots in $S^{n+2}$ and the set of
$\Q$-concordance classes of $n$-knots in rational homology spheres,
respectively. They are abelian groups under connected sum.  We have a
natural homomorphism $C^{\vphantom{\Q}}_n\to C^\Q_n$ since concordant
links in a sphere are $\Q$-concordant.

\begin{thm}\label{thm:independent-knots}
There exist infinitely many $(2q-1)$-knots in rational homology
spheres that are linearly independent in the cokernel of
$C^{\vphantom{\Q}}_{2q-1}\to C^\Q_{2q-1}$.
\end{thm}

We remark that knots constructed by Cochran and Orr using the covering
link construction~\cite{CO2} are nontrivial in the cokernel since they
are not concordant to knots of complexity one.  Our result shows that
the cokernel is large enough to contain a subgroup isomorphic
to~$\Z^\infty$.

On the other hand, we illustrate that our theory is also useful for
links in spheres by offering alternative proofs of a result of
Litherland~\cite{Li} about signatures of satellite knots and a result
of Kawauchi~\cite{Ka} about concordance of split links. In a separated
paper~\cite{CK4}, we will compute signature jump functions of covering
links and study concordance of boundary links and homology boundary
links.  In particular, results of Cochran and Orr~\cite{CO1,CO2},
Gilmer and Livingston~\cite{GL}, and Levine~\cite{L4} will be
generalized.

\section{Signatures of matrices}
\label{sec:matrix-sign-poly}

In this section, we study algebraic properties of the signature
function $\sigma^q_A$ and its jump function $\delta^q_A$ that will
play a key role in our link signature theory.  It is known that if $A$
is a complex square matrix such that $A\pm \bar A^T$ is nonsingular,
the function $e^{i\phi}\mapsto \sigma^q_A(\phi)$ is constant except at
finitely many jumps where $e^{i\phi}$ is a zero of $\det(tA-\epsilon
\bar A^T)$, and so the jump function $\delta^q_A$ is well-defined. In
general, $\delta^q_A$ is well-defined for any complex square matrix
$A$. We accompany a proof since we did not find one in the literature.
Let $\mathstrut^-$ denote the usual involution $(\sum \alpha_i t^i)^-
=\sum \bar{\alpha_i}t^{-i}$ on the Laurent polynomial
ring~$\C[t,t^{-1}]$.

\begin{lem}\label{lem:sign-local-const}
Let $Q(t)$ be a square matrix over $\C[t,t^{-1}]$ that satisfies
$Q(t)=(Q(t)^T)^-$.  Then for any $\phi_0 \in \R$, the signature of
$Q(e^{i\phi})$ is locally constant on a deleted neighborhood of
$\phi_0$, i.e.,\ it is constant on each of $(\phi_0-r,\phi_0)$ and
$(\phi_0,\phi_0+r)$ for some~$r>0$.
\end{lem}

\begin{proof}
First we show that the rank of $Q(e^{i\phi})$ is constant on a deleted
neighborhood of $\phi_0$.  Let $R$ be a square matrix obtained by
deleting some rows and columns from $Q(t)$ and denote its determinant
by $d_R(t)$. Then either $d_R(t)\equiv 0$ or there is a deleted
neighborhood $U_R$ of $\phi_0$ where $d_R(e^{i\phi})\ne 0$ since
$d_R(t)$ is a Laurent polynomial over $\C$.  Since there are finitely
many $R$, the intersection of all $U_R$ is a deleted neighborhood of
$\phi_0$ where $\rank Q(e^{i\phi})$ is constant.

Now it suffices to show that $\sign Q(e^{i\phi})$ is locally constant
on any open interval $J$ where $\rank Q(e^{i\phi})$ is constant.  Near
a fixed $\phi_1$ in $J$, $Q(e^{i\phi})$ is congruent to a diagonal
matrix whose diagonals $\lambda_i(\phi)$ are real-valued continuous
functions on $\phi$. We can choose a neighborhood $U$ of $\phi_1$ in
$J$ such that if $\lambda_j(\phi_1)\ne 0$, then $\lambda_j(\phi)\ne 0$
and its sign is not changed in $U$. Since $\rank Q(e^{i\phi})$ is
constant, $\lambda_j(\phi)=0$ on $U$ if
$\lambda_j(\phi_1)=0$. Therefore $\sign Q(e^{i\phi})$ is constant on
$U$.
\end{proof}

From the lemma, it is immediate that for an arbitrary $A$, the
function $e^{i\phi}\mapsto \sigma^q_A(\phi)$ is locally constant
except at finitely many jumps on the unit circle and so $\delta^q_A$
is well-defined.

As mentioned in the introduction, signature invariants of links are
defined from the signature functions of Seifert matrices of
generalized Seifert surfaces. Note that a union of parallel copies of
a generalized Seifert surface is again a generalized Seifert
surface. Since both of them can be used to define signature invariants
of links, we are naturally led to investigate a relation of the
signature functions of them.  Explicitly, for an $r$-tuple
$\alpha=(s_1,\ldots,s_r)$ with $s_n=\pm 1$, consider the union of $r$
parallel copies of a Seifert surface $F$ of a $(2q-1)$-link, where the
$n$-th copy is oriented according to the sign of $s_n$. If $A$ is a
Seifert matrix of $F$, the Seifert matrix of the union of parallel
copies is given as follows.
$$
i^q_\alpha A=
\begin{bmatrix}
A_1 & A & A & \cdots & A\\ \epsilon A^T & A_2 & A & \cdots & A\\ \epsilon
A^T & \epsilon A^T & A_3 & \cdots & A \\ \vdots & \vdots &  \vdots & \ddots
& \vdots \\ \epsilon A^T & \epsilon A^T & \epsilon A^T & \cdots & A_r
\end{bmatrix}
$$
where
$$
A_n = \left\{\begin{array}{ll}
A & \mbox{ if }s_n = 1,\\
\epsilon A^T & \mbox{ if }s_n = -1.
\end{array}\right.
$$
The following result shows that the signature jump functions of
$i^q_\alpha A$ and $A$ are essentially equivalent up to a change of
parameters. To simplify statements, we define $\sgn 0$ to be zero.

\begin{thm}[Reparametrization formula]\label{thm:sign-formula}
If $n_\alpha\ne 0$ or $A-\epsilon A^T$ is nonsingular, then
$\delta^q_{i^q_\alpha A}(\theta) = \sgn n_\alpha \cdot
\delta^q_A(n_\alpha \theta)$.
\end{thm}

The rest of this section is devoted to the proof of
Theorem~\ref{thm:sign-formula}.

For later use, we introduce some notations. For a nonzero integer $r$,
let $i^q_rA=i^q_{\alpha^r}A$ where $\alpha^r$ is an $|r|$-tuple whose
entries are $r/|r|$.  We denote $i^1_{\alpha}$, $i^2_{\alpha}$ by
$i^+_{\alpha}$, $i^-_{\alpha}$, respectively, and use similar
notations for $\sigma$ and $\delta$.  Since Seifert pairing of links
in rational homology spheres is rational-valued (see
Section~\ref{sec:link-in-rhs}), we assume that $A$ is a square matrix
over~$\Q$.  For a positive integer $a$, let $\eta_a \colon \R \to \Z$
be the step function of period $2\pi$ given by
$$
\eta_a(\phi) = \begin{cases}
a+1-2k, & \text{if }2\pi(k-1)/a <\phi <2\pi k/a \quad (k=1,\ldots,a)\\
a-2k, & \text{if }\phi = 2\pi k/a \quad (k=1,\ldots,a-1)\\
0 & \text{if }\phi=0.
\end{cases}
$$
For negative $a$, let $\eta_{a}(\phi) = \eta_{-a}(-\phi)$ and let
$\eta_0(\phi)=0$.

\begin{lem}\label{lem:alg-sign-formula}
Suppose that $e^{ik\phi}\ne 1$ for $k = 1,\ldots,r$.  If $n_\alpha \ne 0$
or $A-\epsilon A^T$ is nonsingular, we have
\begin{enumerate}
\item $\sigma^+_{i^+_\alpha A}(\phi) = \sigma^+_A(n_\alpha \phi)$.
\item $\sigma^-_{i^-_\alpha A}(\phi) =
\eta_{n_\alpha}(\phi)\cdot \sign(A+A^T) + \sigma^-_A(n_\alpha \phi)$.
\end{enumerate}
\end{lem}

\begin{proof}
Let $w = e^{i\phi}$. First we note that if $w\ne 1$,
$$
\sigma^q_A(\phi)=\sign\Big(\frac{wA-\epsilon \bar
A^T}{w-1}\Big).
$$
We will prove some matrix congruence relations. Consider the matrix
$$
\begin{bmatrix}
\frac{w^kA-\epsilon A^T}{w^k-1} & A & \bar X^T \\
\epsilon A^T & \frac{wB-\epsilon B^T}{w-1} & \bar X^T \\
X & X & Y
\end{bmatrix}
$$
where $B=A$ or $\epsilon A^T$.  Subtracting the second row from the
first row and performing the corresponding column operation, a new
matrix is obtained. Subsequently we perform the following row
operations and corresponding column operations: If $k \ne -1$ and
$B=A$, we add the first row times $w(w^k-1)/(w^{k+1}-1)$ to the second
row.  If $k = -1$ and $B=A$, we assume $A-\epsilon A^T$ is nonsingular
and add the top row times $(1-w)X(A-\epsilon A^T)^{-1}$ to the bottom
row.  Then we obtain
$$
\begin{bmatrix}
\frac{(w^{k+1}-1)(A-\epsilon A^T)}{(w-1)(w^k-1)} & 0 & 0 \\
0 & \frac{w^{k+1}A-\epsilon A^T}{w^{k+1}-1} & \bar X^T \\
0 & X & Y
\end{bmatrix}  \text{ and }
\begin{bmatrix}
0 & \frac{\epsilon A^T-A}{w-1} & 0 \\
\frac{w(\epsilon A^T-A)}{w-1} & \frac{wA-\epsilon A^T}{w-1} & 0 \\
0 & 0 & Y
\end{bmatrix}
$$
respectively.  If $B=\epsilon A^T$, replace $w$ by $w^{-1}$ and $k$ by
$-k$, respectively.  Then we obtain
$$
\begin{bmatrix}
\frac{(w^{k-1}-1)(A-\epsilon A^T)}{(w^{-1}-1)(w^k-1)} & 0 & 0 \\
0 & \frac{w^{k-1}A-\epsilon A^T}{w^{k-1}-1} & \bar X^T \\
0 & X & Y
\end{bmatrix} \text{ and }
\begin{bmatrix}
0 & \frac{w(A-\epsilon A^T)}{w-1} & 0 \\
\frac{A-\epsilon A^T}{w-1} &
\frac{w(\epsilon A^T)-A}{w-1} & 0
\\ 0 & 0 & Y
\end{bmatrix}
$$
for $k\ne1$ and $k=1$, respectively.

Now we prove the lemma.  Note that $\epsilon (i^q_\alpha A)^T =
i^q_{-\alpha} A$ and $\sigma^q_A(-\phi)=\sigma^q_{\epsilon
A^T}(\phi)$, where $-(s_1,\ldots,s_r) = (-s_1,\ldots,-s_r)$.  We may
assume that $n_\alpha \ge 0$ and furthermore by rearranging $s_i$ we
may assume $s_1=\cdots=s_u=1$ and $s_{u+1}=\cdots=s_r=-1$ where $u \ge
r/2$.  For, if $n_\alpha < 0$ and the lemma is proved for the case of
$n_\alpha \ge 0$,
\begin{multline*}
\sigma^q_{i^q_\alpha A}(\phi) =
\sigma^q_{\epsilon(i^q_{(-\alpha)}A)^T}(\phi) =
\sigma^q_{i^q_{-\alpha}A}(-\phi) \\
=\begin{cases}
\sigma^+_A((-n_\alpha)(-\phi)) &
\text{if }\epsilon=1\\
\sigma^-_A((-n_\alpha)(-\phi))+\eta_{-n_\alpha}(-\phi)\sign(A+A^T) &
\text{if }\epsilon=-1\\
\end{cases}
\end{multline*}
and the conclusion follows.

Suppose $n_\alpha>0$. By applying repeatedly the first and third
congruence relations in the above to the matrix
\begin{equation}\label{equ:cablematrix}
\frac{w(i^q_\alpha A)-(i^q_\alpha A)^T}{w-1}
=
\begin{bmatrix}
A_1 & A & A & \cdots & A\\
\epsilon A^T & A_2 & A & \cdots & A\\
\epsilon A^T & \epsilon A^T & A_3 & \cdots & A \\
\vdots & \vdots &  \vdots & \ddots & \vdots \\
\epsilon A^T & \epsilon A^T & \epsilon A^T & \cdots & A_r
\end{bmatrix}
\end{equation}
where
$$
A_n = \left\{\begin{array}{ll}
\frac{wA-\epsilon A^T}{w-1} & \mbox{ if }s_n = 1\\
\frac{w(\epsilon A^T)-\epsilon(\epsilon A^T)^T}{w-1} & \mbox{ if }s_n = -1
\end{array}\right.,
$$
we obtain a matrix that is a block sum of
\begin{multline}\label{equ:reducedmatrix}
\frac{(w^2-1)(A-\epsilon A^T)}{(w-1)(w-1)},
\ldots,
\frac{(w^u-1)(A-\epsilon A^T)}{(w-1)(w^{u-1}-1)},\\
\frac{(w^{u-1}-1)(A-\epsilon A^T)}{(w^{-1}-1)(w^u-1)},
\ldots,
\frac{(w^{u-(r-u)}-1)(A-\epsilon A^T)}{(w^{-1}-1)(w^{u-(r-u)+1}-1)}
\end{multline}
and
$$
\frac{w^{n_\alpha}A-\epsilon A}{w^{n_\alpha}-1}.
$$

If $\epsilon=1$, the contribution of (\ref{equ:reducedmatrix}) to the
signature is zero since $\sign(i(A-A^T))=0$. (Recall that $A$ is a
matrix over $\Q$.) This proves the conclusion for $n_\alpha>0$ and
$\epsilon=1$.

If $\epsilon=-1$, let $f_n(\phi) = -\sin(n+1)\phi + \sin n\phi + \sin
\phi$. Then we have
$$
\operatorname{sgn}f_n(\phi)=
\operatorname{sgn} \frac{i(w^{n+1}-1)}{(w-1)(w^n-1)} =
-\operatorname{sgn} \frac{i(w^{n}-1)}{(w^{-1}-1)(w^{n+1}-1)}
$$
and therefore
\begin{align*}
\sigma^q_{i^q_A}(\phi) &=
 \sigma^q_A(n_\alpha \phi) +
 \sign(A+A^T) \Big(\sum_{n=1}^{u-1} \operatorname{sgn} f_{n}(\phi) -
 \sum_{n=1}^{r-u} \operatorname{sgn} f_{u-n}(\phi) \Big)
 \\
&=\sigma^q_A(n_\alpha \phi) +
 \sign(A+A^T) \sum_{n=1}^{n_\alpha-1} \operatorname{sgn} f_{n}(\phi).
\end{align*}
By an elementary calculation, it is easily shown that $f_n$ is zero at
$\phi = 2\pi k/(n+1)$, $2\pi k/n$, and $f_n(\phi)>0$ on $(2\pi
(k-1)/n, 2\pi k/(n+1))$ and $f_n(\phi)<0$ on $(2\pi k/(n+1), 2\pi
k/n)$. Hence by induction on $n_\alpha$ the sum in the above equation
is equal to $\eta_{n_\alpha}(\phi)$. This proves the conclusion for
$n_\alpha>0$ and $\epsilon=-1$.

Suppose $n_\alpha=0$. In this case, the proof also proceeds in a
similar way. We transform the matrix (\ref{equ:cablematrix}) to a
block sum of matrices in (\ref{equ:reducedmatrix}) and
\begin{equation}\label{equ:hyperbolicmatrix}
\begin{bmatrix}
0 & \frac{w}{w-1}(A-\epsilon A^T) \\
\frac{1}{w-1}(A-\epsilon A^T) & \frac{w(\epsilon A^T)-A}{w-1}
\end{bmatrix}
\end{equation}
by applying all of the four congruence relations shown before.  The
signature of (\ref{equ:hyperbolicmatrix}) is zero since $A-\epsilon
A^T$ is nonsingular. Hence we have $\sigma^q_{i^q_\alpha A}(\phi) =
0$.  On the other hand, $\sigma^+_A(n_\alpha\phi) = \sigma^+_A(0) = 0$
for $\epsilon=1$, and
$\eta_{n_\alpha}(\phi)\sign(A+A^T)+\sigma^-_A(n_\alpha \phi) = 0$ for
$\epsilon=-1$. This completes the proof.
\end{proof}

Now we are ready to prove the main theorem of this section.

\begin{proof}[Proof of Theorem~\ref{thm:sign-formula}]
If $n_\alpha=0$ and $A-\epsilon A^T$ is nonsingular,
$\sigma^q_{i^q_\alpha A}$ is the zero function by
Lemma~\ref{lem:alg-sign-formula}, and so $\delta^q_{i^q_\alpha
A}(\theta) = 0 = \sgn n_\alpha \cdot \delta^q_A(n_\alpha \theta)$ for
all $\theta$.

Suppose $n_\alpha>0$.  By Lemma~\ref{lem:sign-local-const}, we can
choose $\epsilon_0, \epsilon_1 > 0$ such that $\sigma^q_\phi(A)$ is
constant on each of $[\theta-\epsilon_0, \theta)$, $(\theta,
\theta+\epsilon_1]$, $[n_\alpha \theta-n_\alpha \epsilon_0, n_\alpha
\theta)$ and $(n_\alpha \theta, n_\alpha \theta+n_\alpha
\epsilon_1]$. We may assume that $(\theta-\epsilon_0)/\pi$ and
$(\theta+\epsilon_1)/\pi$ are irrational.  For $\epsilon=1$, by
Lemma~\ref{lem:alg-sign-formula}
\begin{align*}
\delta^+_{i^+_\alpha A}(\theta) &=
\sigma^+_{i^+_\alpha A}(\theta+\epsilon_1)
-\sigma^+_{i^+_\alpha A}(\theta-\epsilon_0) \\
&= \sigma^+_A(n_\alpha (\theta+\epsilon_1))
-\sigma^+_A(n_\alpha (\theta-\epsilon_0)) \\
&= \delta^+_A(n_\alpha \theta).
\end{align*}
For $\epsilon=-1$, if $\theta=2\pi k$ then the conclusion is trivial
since $\delta^-_{A}(2\pi k)=0$ for any $A$. If not, we have
\begin{equation*} \begin{split}
\delta^-_{i^-_\alpha A}(\theta) &=
\big(\eta_{n_\alpha}(\theta+\epsilon_1)-\eta_{n_\alpha}(\theta-\epsilon_0)\big)
\sign(A+A^T)\\
&\qquad\qquad + \big(\sigma^-_{A}(n_\alpha \theta+n_\alpha \epsilon_1)
-\sigma^-_{A}(n_\alpha \theta-n_\alpha \epsilon_0) \big).
\end{split}\end{equation*}
by Lemma~\ref{lem:alg-sign-formula}. If $\theta \ne 2\pi k/n_\alpha$
for any $k$, we may assume that $\eta_{n_\alpha}$ is constant near
$\theta$, and the theorem is proved from the above equation.  If
$\theta=2\pi k/n_\alpha$ for some $k$ that is not a multiple of
$n_\alpha$, then $\delta^-_{i^-_\alpha A}(\theta) = -2\sign(A+A^T)+
\big(\sigma^-_A(2\pi k + n_\alpha \epsilon_1) -\sigma^-_A(2\pi k -
n_\alpha \epsilon_0) \big) = 0 = \delta^-_A(n_\alpha\theta)$ since
$\sigma^-_A(\pm\phi) = \pm\sign(A+A^T)$ for all sufficiently small
$\phi>0$.

If $n_\alpha < 0$, we have $\delta^q_{i^q_\alpha}(\theta) =
\delta^q_{\epsilon (i^q_{-\alpha}A)^T}(\theta) =
-\delta^q_{i^q_{-\alpha}A}(-\theta) = -\delta_A^q(n_\alpha\theta)$, as
in the proof of Lemma~\ref{lem:alg-sign-formula}.
\end{proof}

We remark that a similar argument shows an analogous result for
Alexander polynomials defined by $\Delta^q_A(t) = \det(tA-\epsilon
A^T)$: If $A-\epsilon A^T$ is nonsingular, then $\Delta^q_{i^q_\alpha
A}(t) = \Delta^q_A(t^{n_\alpha})$ up to units of $\Q[t,t^{-1}]$. This
result is mentioned and used in~\cite{CK1}.

We finish this section with a remark on the matrix cobordism groups.
Consider the set of rational matrices $A$ such that $A-\epsilon A^T$
is nonsingular.  Define the block sum operation and cobordisms of such
matrices as in~\cite{L2}. Then arguments of~\cite{L2} show that the
cobordism of such matrices is an equivalence relation and that the set
$G^{\Q}_\epsilon$ of cobordism classes becomes an abelian group under
the block sum.  Then $i^q_\alpha$ induces a group homomorphism on
$G^\Q_\epsilon$ and both $\delta^q_A$ and $i^q_\alpha$ are
well-defined on $G^\Q_\epsilon$. Thus Lemma~\ref{lem:alg-sign-formula}
also makes sense when $A$ is regarded as its cobordism class.

\section{Seifert surfaces of links}
\label{sec:link-in-rhs}

In this section, we study links and their Seifert surfaces in rational
homology spheres. We start with a review of rational-valued linking
numbers.  Let $x$, $y$ be disjoint $p$, $q$-cycles respectively in a
$(n+2)$-rational homology sphere $\Sigma$ where $p+q=n+1$. Since
$H_p(\Sigma)$, $H_q(\Sigma)$ are torsion, there are $(p+1)$,
$(q+1)$-chains $u$, $v$ in $\Sigma$ such that $\partial u = ax,
\partial v = by$ for some nonzero integers $a$, $b$.  The
\emph{linking number} of $x$ and $y$ is defined to be $\lk_\Sigma(x,
y) = (1/b)(x\cdot v) = (-1)^{q+1}(1/a)(u\cdot y)$ where~$\cdot$
denotes the intersection number in $\Sigma$. The linking number
$\lk_\Sigma(x,y)$ is determined by $x$ and $y$ and is independent of
the choices of $a$, $b$, $u$ and $v$.  When the ambient space is
obvious, we will denote the linking number by~$\lk(x,y)$.  We remark
that $\lk(-,-)$\ is well-defined only for disjoint cycles, and its
value modulo $\Z$ is well-defined for homology classes.

The following properties of the linking number are easily
checked. First, if $c(x'-x)=\partial w$ for some $c\ne 0$ and
$(p+1)$-chain $w$ disjoint to $y$, then $\lk(-,y)$ assumes the same
value on $x, x'$.  $\lk(x,-)$ has the same property.  Second, suppose
that $\Sigma$ is a boundary component of an $(n+3)$-manifold $\Delta$
such that $H_{p+1}(\Delta;\Q)=H_{q+1}(\Delta;\Q)=0$.  If $x$, $y$ are
as above and $u$, $v$ are relative $(p+1)$, $(q+1)$-cycles of
$(\Delta, \Sigma)$ such that $\partial u = ax, \partial v = by$ for
some $a,b\ne 0$, $\lk(x,y)=(1/ab)(u\cdot v)$.

We give a formula for linking numbers in a rational homology 3-sphere
described by a Kirby diagram. For an $m$-component framed link, we
define the linking matrix $A$ to be a square matrix of dimension $m$
whose $(i,j)$-entry is the linking number of the $i$-th component and
the $j$-th longitude with respect to the framing.

\begin{thm}\label{thm:lk-kirby-diagram-dim1}
Let $L=K_1\cup\cdots\cup K_m$ be a framed link in $S^3$ such that the
result of surgery on $S^3$ along $L$ is a rational homology
sphere~$\Sigma$. Let $A=(a_{ij})$ be the linking matrix of $L$.  For
disjoint $1$-cycles $a$, $b$ in $S^3-L$, we have
$$
\lk_\Sigma(a,b)=\lk_{S^3}(a,b)-x^T A^{-1}y
$$
where $x=(x_i)$, $y=(y_i)$ are column vectors with $x_i=\lk_{S^3}(a,
K_i)$, $y_i=\lk_{S^3}(b, K_i)$, respectively.
\end{thm}
\begin{proof}
Let $\mu_i$ be a meridian curve on a tubular neighborhood of
$K_i$. Let $u$ and $u_i$ be 2-chains in $S^3$ that are transverse to
$L$ and bounded by $b$ and a 0-linking longitude of $K_i$,
respectively. Let $\bar u$ and $\bar u_i$ be chains obtained by
removing from $u$ and $u_i$ the intersection with an open tubular
neighborhood of $L$, respectively. Let $c_i$ be the core of the
2-handle attached along $K_i$.  Let $v=k\bar u+\sum_i
z_i(\bar{u_i}+c_i)$ be a chain in $\Sigma$, where $k$, $z_i$ are
integers to be specified later. Then $\bd v$ is homologous to
$kb-\sum_j \sum_i(ky_j+z_i a_{ij})\mu_j$. Thus if $z=(z_1,\ldots,z_m)$
satisfies $z A=-ky$, then $\bd v$ is homologous to $kb$. Since
$\Sigma$ is a $\Q$-homology sphere, $A$ is invertible over $\Q$. So we
can take as $k$ a nonzero common multiple of denominators of entries
of $A^{-1}y$, and $z=-k A^{-1}y$. Now we have
\begin{equation*}
\begin{split}
\lk_\Sigma(a,b) &= \frac{1}{k}(a\cdot v)\\
&= \frac{1}{k}\Big(\lk_{S^3}(a,kb)+\sum_i x_i z_i\Big) \\
&= \lk_{S^3}(a,b)-x^T A^{-1} y
\end{split}
\end{equation*}
\end{proof}

Now we turn our attention to links.  The main goal of this section is
to study the existence of a generalized Seifert surface of a given
type $\tau$ for links. Let $L$ be an $n$-link in a rational homology
sphere, and suppose $n>1$. $L$ can be viewed as a framed manifold
since all framings are equivalent. Generally, for a framed submanifold
$M$ of codimension~2 in a manifold $W$, we identify a tubular
neighborhood of $M$ with $M\times D^2$ with respect to the framing,
and in particular identify $M\times S^1$ with a submanifold in the
boundary of the exterior~$E_M = W-M\times \inte(D^2)$.  We say $(W,M)$
is \emph{primitive} if there is a homomorphism of $H_1(E_M) \to \Z$
whose composition with $H_1(M\times S^1) \to H_1(E_M)$ is the same as
the projection $H_1(M\times S^1)\to H_1(S^1)=\Z$. When $W$ is obvious,
we say $M$ is primitive. By a standard use of obstruction theory,
transversality and Thom-Pontryagin construction, $L$ is primitive if
and only if $L$ bounds a Seifert surface.

Our main observation is that certain unions of parallel copies of $L$
is primitive and hence bound Seifert surfaces. We use the following
notations to denote parallel copies of a framed submanifold $M$
in~$W$.  For an $r$-tuple $\alpha=(s_1,\ldots,s_r)$ where $s_i=\pm1$,
let $n_\alpha$ be the sum of $s_i$ as before, and $i_\alpha M$ be the
union of $r$ parallel copies of $M$ in $W$ with respect to the framing
where the $i$-th copy is oriented according to the sign of $s_i$.  For
a nonzero integer $r$, let $i_rM=i_\alpha M$ where $\alpha$ is the
$r$-tuple $(\sgn r,\ldots,\sgn r)$.

We assert that $i_\alpha L$ is primitive if and only if $n_\alpha$ is
a multiple of~$c$, where $c$ is a positive integer given as follows.
Let $\varphi\colon H_1(E_L;\Q) \to \Q$ be the homomorphism sending
each meridian to~1, which is uniquely determined by the Alexander
duality. Let $P$ be the quotient group of $H_1(E_L)$ by its torsion
subgroup. Since $P\otimes \Q = H_1(E_L;\Q)$, we can view $P$ as a
subgroup of $H_1(E_L;\Q)$. Choose a basis $\{e_i\}$ of $P$, and let
$c$ be the positive least common multiple of denominators of the
reduced fractional expression of $\varphi(e_i)$.  We call $c$ be the
\emph{complexity} of~$L$. $c$ is independent to the choice of the
basis.  Note that $c\varphi(P) \subset \Z \subset \Q$ and hence a
homomorphism $\phi$ from $P$ into $\Z$ is induced by~$c\varphi$.

Suppose that $i_\alpha L$ is primitive for $\alpha=(s_1,\ldots,s_r)$.
Assuming that $i_\alpha L$ is the product of $L$ and $r$ points in the
tubular neighborhood $L\times D^2$, we can view $E_L$ as a subspace of
$E_{i_\alpha L}$. There is a homomorphism of $H_1(E_{i_\alpha L})$
into $\Z$ sending each meridian of $i_\alpha L$ to~1, and it induces a
homomorphism $h$ of $P=H_1(E_L)/\text{torsion}$ to $\Z$ sending each
meridian of $L$ to $n_\alpha$. The homomorphisms
$(1/n_\alpha)(h\otimes 1_{\Q})$ and $\varphi$ on $P\otimes
\Q=H_1(E_L;\Q)$ are equal since they assume the same value on
meridians.  Hence $h(e_i)/n_\alpha = \varphi(e_i)$ and so $n_\alpha$
is a multiple of the denominator of $\varphi(e_i)$. This shows that
$n_\alpha$ is a multiple of~$c$. Conversely, suppose $n_\alpha = kc$
for some integer $k$. Then the homomorphism $k\phi\colon H_1(E_L) \to
\Z$ sends each meridian of $L$ to $n_\alpha$, and it extends to a
homomorphism of $H_1(E_{i_\alpha L})$ to $\Z$ sending each meridian of
$i_\alpha L$ to~1. This proves the assertion.

For a type $\tau=(\tau_i)$ for $L$, let $L^\tau$ be the union of
$\tau_i$ parallel copies (negatively oriented if $\tau_i<0$) of the
$i$-th component of $L$ over all~$i$.  We call the complexity of
$L^\tau$ the \emph{complexity of type $\tau$ for $L$} and denote it by
$c^\tau(L)$.  A generalized Seifert surface of type $\tau$ with
complexity $r$ for $L$ is a Seifert surface of $i_r L^\tau$, which
exists if and only if $i_rL^\tau$ is primitive. This holds exactly
when $r$ is a multiple of $c^\tau(L)$. In particular,

\begin{thm}\label{thm:seifert-surface-existence-highdim}
For $n>1$, an $n$-link $L$ admits a generalized Seifert surface of any
type~$\tau$.
\end{thm}

\begin{thm}\label{thm:compute-complexity}
Fix a basis of $P=H_1(E_L)/\text{torsion}$, and let $B$ be a matrix
whose $i$-th column represents the $i$-th meridian of $L$ with respect
to the basis. Then $c^\tau(L)$ is equal to the least common multiple
of denominators (of a reduced fractional expression) of entries of the
row vector $\tau^T B^{-1}$. (We view $\tau$ as a column vector.)
\end{thm}

\begin{proof}
First of all, $B$ is invertible since meridians are independent.
Consider the homomorphism of $H_1(E_{L^\tau};\Q)$ to $\Q$ sending
meridians to~1. Viewing $E_L$ as a subspace of $E_{L^\tau}$, a
homomorphism of $H_1(E_L;\Q)$ to $\Q$ sending the $i$-th meridian of
$L$ to $\tau_i$ is induced. Suppose that the $i$-th basis of $P$ is
sent to $x_i$ under the induced homomorphism. By the definition of the
complexity, the least common multiple of the denominators of $x_i$ is
equal to $c^\tau(L)$.  The row vector $x=(x_i)$ satisfies $xB=\tau^T$,
and so $x=\tau^T B^{-1}$ as desired.
\end{proof}

Now we consider 1-links in a rational homology sphere. Our result for
1-links is expressed in terms of a linking matrix of $L$ (for any
choice of a framing) as follows.  We say a 1-link \emph{$L$ admits
type $\tau$} if there is a generalized Seifert surface of type $\tau$
for~$L$.

\begin{thm}\label{thm:seifert-surface-existence-dim1}
A 1-link $L$ admits type $\tau$ if and only if for a linking matrix
$A=(a_{ij})$ with respect to a (or any) framing,
$$
\frac{1}{\tau_i} \sum_j a_{ij}\tau_j
$$
is an integer for all $i$.
\end{thm}

First we will investigate when a generalized Seifert surface inducing
a given fixed framing exists, and then show
Theorem~\ref{thm:seifert-surface-existence-dim1}.  Let $L$ be a
1-link, $f$ be a framing on $L$ and $A=(a_{ij})$ be the linking matrix
with respect to~$f$.

\begin{lem}\label{lem:seifert-surface-existence-framed}
There exists a generalized Seifert surface of type $\tau=(\tau_i)$ for
$L$ that induces the same framing on its boundary as the framing
induced by $f$ if and only if
$$
\sum_{j} a_{ij}\tau_j = 0
$$
for all $i$.
\end{lem}

\begin{proof}
Let $L^\tau$ be the link obtained by taking $\tau_i$ parallel copies
of the $i$-th components of $L$ with respect to the given framing~$f$.
$L$ admits a desired generalized Seifert surface of complexity $r$ if
and only if $i_rL^\tau$ is primitive.  Since $H_1(E_{L^\tau};\Q)$ is
generated by meridians of $L^\tau$, for any multiple $r$ of the
complexity of $L^\tau$ we obtain a homomorphism $\phi$ of $H_1(E_{i_r
L^\tau})$ to $\Z$ sending each meridian of $i_r L^\tau$ to 1, by the
arguments for higher dimensional links.  In order to show that $i_r
L^\tau$ is primitive, we need to check whether the longitudes of $i_r
L^\tau$ is sent to zero by $\phi$, or equivalently, the longitudes of
$L$ is sent to zero by the unique homomorphism from $H_1(E_L;\Q) \to
\Q$ sending meridians $\mu_j$ of $L$ to~$\tau_j$. Since the $i$-th
longitude of $L$ is equal to $\sum_j a_{ij} \mu_j$ in $H_1(E_L;\Q)$,
the conclusion follows.
\end{proof}

\begin{proof}[Proof of Theorem~\ref{thm:seifert-surface-existence-dim1}]
Let $A$ be a linking matrix of a 1-link $L$ with respect to a fixed
framing~$f$.  Then the $i$-th longitude with respect to an arbitrary
framing $f'$ is homologous to the sum of the $i$-th longitude with
respect to $f$ and a multiple of the $i$-th meridian in the complement
of $L$.  Hence the linking matrix with respect to $f'$ is of the form
$A+D$ where $D$ is a diagonal matrix of integers.  Conversely, $D$
determines a framing since framings are in 1-1 correspondence with
elements of $H_1(L;\pi_1(SO_2))\cong \Z^m$.

Therefore, by Lemma~\ref{lem:seifert-surface-existence-framed}, there
is a generalized Seifert surface of type $\tau$ if and only if
$(A+D)\tau=0$ for some diagonal matrix $D$ of integers, or
equivalently, $\frac{1}{\tau_i}\sum_j a_{ij}\tau_j$ is an integer for
all~$i$.
\end{proof}

We remark that generalized Seifert surfaces of the same type induce
the same framing, since in the above proof the numbers $d_i$ are
determined by $\tau=(\tau_i)$. From now on, for a 1-link $L$ and a
type $\tau$ admitted by $L$, we denote by $L^\tau$ the union of
$\tau_i$ parallel copies of the $i$-th component of $L$ with respect
the framing determined by $\tau$.

\begin{exmp}
1. If a linking matrix of a link $L$ is zero, then any type is admitted.

2. If $L$ is a link in a homology sphere, each entry of a linking
matrix is an integer. Hence type~$(1\cdots 1)$ is admitted.

3. Let $C$ be an $m$-component unlink in $S^3$ and $L$ be the union of
meridian curves of components of $C$.  A rational homology sphere
$\Sigma$ is obtained by $n$-surgery on each component of $C$ for $n\ge
2$, and $L$ can be viewed as a link in $\Sigma$.  By
Theorem~\ref{thm:lk-kirby-diagram-dim1}, the diagonal matrix with
diagonals $-1/n$ is a linking matrix, and by
Theorem~\ref{thm:seifert-surface-existence-dim1}, $L$ does not admit
any type.

4.  For any given type $\tau=(a,b)$, there is a link that admits only
type $\pm\tau$. For, suppose that a link $L$ has linking matrix
$\spmatrix{-b/3a & 1/3 \\ 1/3 & -a/3b}$. If a type $(x,y)$ is admitted
by $L$, $(y/x-b/a)/3$ and $(x/y-a/b)/3$ are integers by
Theorem~\ref{thm:seifert-surface-existence-dim1}. This implies that
$(x,y)=\pm (a,b)$. By the following lemma that classifies linking
matrices, such a link $L$ exists.
\end{exmp}

\begin{lem}
For any symmetric square matrix $A$ with rational entries, there is a
framed 1-link in a rational homology sphere whose linking matrix
is~$A$.
\end{lem}

\begin{proof}
We claim that for any symmetric nonsingular $m$ by $m$ matrix $V$ and
$m$ by $n$ matrix $B=(b_{ij})$ with integral entries, there is an
$n$-component framed 1-link in a rational homology sphere with linking
matrix $B^T V^{-1} B$.  Consider the $n$-component standard unlink $L$
in $S^3$ with zero-linking framing. We can construct an $m$-component
framed link $L'$ in $S^3$ with linking matrix $-V$ such that $L$ and
$L'$ are disjoint and the linking number of the $i$-th component $L'$
and the $j$-th component of $L$ is $b_{ij}$. The result of surgery
along $L'$ is a rational homology 3-sphere $\Sigma$ since $-V$ is
nonsingular. $L$ can be viewed as a link in $\Sigma$, and the linking
matrix of $L$ in $\Sigma$ is $B^T V^{-1} B$, by
Theorem~\ref{thm:lk-kirby-diagram-dim1}.

Any symmetric matrix $A$ over $\Q$ can be written as $A=P^T
\spmatrix{D & 0 \\ 0 & 0} P$ where $D$ is a nonsingular diagonal
matrix and each $0$ represents a zero matrix of suitable size.  Choose
a nonzero common multiple $n$ of denominators of entries of $P$ and
$D^{-1}$ and let $B=n \spmatrix{I & 0} P$, $V=n^2D^{-1}$. Then $B$ and
$V$ are matrices over $\Z$, and we have $A=B^T V^{-1} B$. By the
claim, $A$ is realized as a linking matrix of a 1-link in a rational
homology sphere.
\end{proof}

We finish our discussion on 1-links by showing that the set of
admitted types of a link is invariant under link concordance. For two
$\Q$-concordant 1-links, a framing on a concordance between the links
induces framings on the links, and in fact this is a 1-1
correspondence between the framings of the two links.

\begin{thm}
If 1-links $L_0$ and $L_1$ are $\Q$-concordant, then $L_0$ admits type
$\tau$ if and only if so does~$L_1$. Furthermore, $\tau$ determines
corresponding framings on $L_0$ and $L_1$.
\end{thm}

\begin{proof}
Since linking numbers are $\Q$-concordance invariants, linking
matrices of $L_0$ and $L_1$ with respect to corresponding framings are
the same. Hence by Theorem~\ref{thm:seifert-surface-existence-dim1},
the conclusion is proved.
\end{proof}

We remark that our complexity is motivated from a similar notion
in~\cite{CO1,CO2}.  The complexity in~\cite{CO1,CO2} is similar to
$c^\tau(L)$ for $\tau=(1\cdots 1)$. In general the former is a
multiple of the latter.

From now on, we will consider only odd dimensional links.  Let
$(\Sigma, L)$ be a $(2q-1)$-link and $F$ be a generalized Seifert
surface of type $\tau$.  Define the \emph{Seifert pairing} $S\colon
H_q(F;\Q)\times H_q(F;\Q) \to \Q$ by $S([x], [y])=\lk(x^+, y)$ for two
$q$-cycles $x$, $y$ in $F$ where $x^+$ denotes the cycle obtained by
pushing $x$ slightly along the positive normal direction of $F$. It is
easily verified that $S$ is well-defined.  Choosing a basis of
$H_q(F;\Q)$, a matrix $A$ over $\Q$ is associated to $S$. $A$ is
called a \emph{Seifert matrix} of $F$. We remark that for $q=1$,
another convention is found in some literature where Seifert pairings
and Seifert matrices are defined on the cokernel of $H_1(\partial
F)\to H_1(F)$. These different conventions give different results in
the case of links (e.g.\ consider a full-twisted annulus with Hopf
link boundary in $S^3$).

We remark that if $A$ is a Seifert matrix of $F$, the matrix
$i^q_\alpha A$ is a Seifert matrix of $i_\alpha F$.

For $q>1$, we can view a Seifert matrix $A$ as a representative of an
element $[A]$ of the matrix cobordism group $G^{\Q}_\epsilon$, since
$A-\epsilon A^T$ represents the intersection pairing on $H_q(F;\Q)$
that is non-degenerated.  The following is an analogous result of the
well-known result of~\cite{L0,L1}.  Recall that a $\Q$-concordance
between two links is called primitive if there is a homomorphism of
the first homology of its exterior into $\Z$ that sends each meridian
and each longitude to~$1$ and~$0$, respectively.

\begin{lem}\label{lem:prim-conc-matrix}
Let $A_0$ and $A_1$ be Seifert matrices of Seifert surfaces for two
$(2q-1)$-links $(\Sigma_0,L_0)$ and $(\Sigma_1,L_1)$ that are
$\Q$-concordant via a primitive $\Q$-concordance. If $q>1$, then
$[A_0] = [A_1]$ in $G^{\Q}_\epsilon$.
\end{lem}

In this case, everything is primitive and an argument in~\cite{L1} can
be used to prove the conclusion.  Details are as follows.

\begin{proof}
Let $(\S, \L)$ be a primitive $\Q$-concordance between
$(\Sigma_0,L_0)$ and $(\Sigma_1,L_1)$, and $E_{\L}$ be the exterior of
$\L$ in $\S$. We may assume $\Sigma_0 \cap E_{\L} = E_{L_0}$ and
$\Sigma_1 \cap E_{\L} = E_{L_1}$.  For $i=0,1$, let $F_i$ be a Seifert
surface of $L_i$ where $A_i$ is defined.  Let $C$ be a parallel copy
of $\L$ in $\bd E_\L$ bounded by $L_1-L_0$.  Let $F = F_0\cup C\cup
-F_1$, a closed oriented submanifold of $\partial E_\L$ of codimension~1.

Define a map $h\colon \partial E_\L \to S^1$ by a Thom-Pontryagin
construction for $F \subset \partial E_\L$. Since $\L$ is primitive,
there is a homomorphism $\phi\colon H_1(E_\L)\to \Z$ sending each
meridian to~1. Let $i_*\colon H_1(\partial E_\L) \to H_1(E_\L)$ denote
the map induced by the inclusion.  Then $\phi i_*$ and $h_*$ assume
the same value on each meridian and so $\varphi i_* = h_*$ since $\Z$
is torsion free.  Therefore $h$ is extended to a map $E_\L\to S^1$. By
a transversality argument, we construct a $(2q+1)$-submanifold $R$ of
$E_\L$ such that $\partial R = F$.

Since $q>1$, $H_q(F) \cong H_q(F_0)\oplus H_q(F_1)$. Let $S$ be the
bilinear pairing on $H_q(F)$ represented by the block sum of $A_0$ and
$A_1$. $H = \Ker(H_q(F;\Q)\to H_q(R;\Q))$ is a half-dimensional
subspace of $H_q(F;\Q)$ by a duality argument, and $S$ vanishes on $H$
by the previous remarks on the linking number as
in~\cite{L1}. Therefore $A_0$ and $A_1$ are cobordant matrices.
\end{proof}

We note that the above proof does not work in the case of $q=1$, since
$H_q(F)$ is not isomorphic to $H_q(F_0)\oplus H_q(F_1)$.

\begin{thm}\label{thm:conc-matrix}
Let $A_1$, $A_2$ be Seifert matrices of two $\Q$-concordant
$(2q-1)$-links with $q>1$. Then $i^q_r[A_0] = i^q_r[A_1]$ in
$G^\Q_\epsilon$ for some integer $r>0$.
\end{thm}

\begin{proof}
Let $\L$ be a $\Q$-concordance between two links whose Seifert
matrices are $A_1$ and $A_2$.  By the same argument as that for links,
there is a nonzero integer $r$ (the ``complexity'' of $\L$) such that
$i_r\L$ is a primitive. $i_r\L$ is a $\Q$-concordance of two links
with Seifert matrices $i^q_rA_0$ and $i^q_rA_1$. By
Lemma~\ref{lem:prim-conc-matrix}, $[i^q_rA_0] = [i^q_rA_1]$.
\end{proof}

\section{Signatures of links}\label{sec:sign-of-link}

In this section we prove the invariance of signatures under link
concordance.  In the case of $q>1$, Seifert matrices can be viewed as
an element of $G^\Q_\epsilon$ and the invariance of signatures may be
shown directly using Theorem~\ref{thm:conc-matrix}.  Because this does
not work for $q=1$, we will interpret signatures via finite branched
covers following ideas of~\cite{K,KT} to show the invariance of
signatures in any odd dimensions.

Suppose that $d$ is a positive integer, $(W,M)$ is a primitive pair,
and $H_1(W;\Q)=H_2(W;\Q)=0$. Then a homomorphism $\phi\colon H_1(E_M)
\to \Z$ sending each meridian to 1 exists uniquely, and the
composition of $\phi$ and the canonical projection $\Z\to\Z_d$
determines a $d$-fold cover $\tilde W$ of $W$ branched along $M$ with
a transformation $t$ on $\tilde W$ corresponding to $1\in\Z_d$ such
that $\tilde W/t\cong W$.  When $W$ is even dimensional, define
$\sigma_{k,d}(W,M)$ to be the signature of the restriction of the
intersection form of $\tilde W$ on the $e^{2k\pi i/d}$-eigenspace of
the induced homomorphism $t_*$ on the middle dimensional homology of
$\tilde W$ with complex coefficients.

Let $L$ be a primitive $(2q-1)$-link (that has the framing induced by
a Seifert surface if $q=1$) in a rational homology sphere $\Sigma$.

\begin{lem}\label{lem:prim-bound}
There is a primitive $(2q+2,2q)$-manifold pair $(W,M)$ such that
$(\Sigma,L)$ is a boundary component of $(W,M)$ and $H_i(W;\Q) =
H_i(\partial W;\Q) = 0$ for $i=1,\ldots,2q$.
\end{lem}

\begin{proof}
Put $W=\Sigma\times[0,1]$ and identify $\Sigma$ with $\Sigma\times 1$.
The map $E_L\to S^1$ given by primitiveness is extended to a map
$\Sigma\to D^2$ in an obvious way, and the union of it and a constant
map $\Sigma\times 0 \to S^1$ is extended to a map $W\to D^2$. Then $M$
is obtained by transversality at an interior point of $D^2$.
\end{proof}

Note that the above argument is similar to that of the well-known
special case $\Sigma=S^{2q+1}$ where $B^{2q+2}$ can be taken
as~$W$. Hence it seems to be natural to take rational homology balls
as $W$ in our general case, however it is not known even when $\Sigma$
is null-cobordant. This is why we allow boundary components of $W$
other than~$\Sigma$.

The following result shows that if $(W,M)$ satisfies the conditions of
Lemma~\ref{lem:prim-bound}, $\sigma_{k,d}(W,M)$ is independent of the
choice of $W$ and $M$.

\begin{lem}\label{lem:sign-welldefined-dim1}
Let $(W_0,M_0)$ and $(W_1,M_1)$ be manifold pairs as in
Lemma~\ref{lem:prim-bound}. Then $\sigma_{k,d}(W_0, M_0) =
\sigma_{k,d}(W_1, M_1)$.
\end{lem}
\begin{proof}
We consider $(W,M)=-(W_0,M_0)\cup_{(\Sigma,L)}(W_1,M_1)$.  The idea is
that $(W,M)$ bounds a primitive pair $(V,N)$ so that the signature of
the branched cover of $W$ vanishes. Then by the additivity of
signatures, the signatures of branched covers of $W_0$ and $W_1$ are
the same and the proof is completed. This is similar to the well-known
argument for the case of $\Sigma=S^{2q+1}$ and $W_i=B^{2q+2}$. In that
case, we can take $V=B^{2q+3}$ and construct $N$ by the argument of
the proof of the previous lemma. For our general case, we take
$V=W\times[0,1]$ and proceed in the same way. Even though $V$ is not
closed, all extra boundary components have vanishing middle
dimensional rational homology and hence have no contributions to
signatures.
\end{proof}

We define $\sigma_{k,d}(L) = \sigma_{k,d}(W,M)$. By
Lemma~\ref{lem:sign-welldefined-dim1}, $\sigma_{k,d}(L)$ is
well-defined.

\begin{lem}\label{lem:sign-branch-cover}
If $L$ is primitive and $A$ is a Seifert matrix of a Seifert surface
of $L$, $\sigma_{k,d}(L) = \sigma^q_A(2 \pi k/d)$.
\end{lem}

In \cite[Chapter XII]{K}, Kauffman shows this lemma when $L$ is a link
in $S^3$ and $\sigma_{k,d}(L)$ is defined from $(D^4,M)$ where $M$ is
a properly embedded surface in $D^4$ bounded by $L$.  In fact the same
argument can be applied to show that if $W=\Sigma\times I$ and $M$ is
obtained by pushing into $W$ the interior of a Seifert surface where
$A$ is defined, then $\sigma_{k,d}(W,M) = \sigma^q_A(2\pi k/d)$. We
omit the details.

The following result as well as Lemma~\ref{lem:prim-conc-matrix} is
based on the idea that most of the techniques of signature invariants
for links in spheres also work for links in rational homology spheres
provided everything is primitive. In fact, for $q>1$, this result is a
consequence of Lemma~\ref{lem:prim-conc-matrix}. We use a geometric
argument to prove this result in any odd dimension.

\begin{lem}\label{lem:prim-conc-sign}
If $A_0$ and $A_1$ are Seifert matrices defined on Seifert surfaces of
primitive $(2q-1)$-links $L_0$ and $L_1$ that are $\Q$-concordant via
a primitive $\Q$-concordance, $\delta^q_{A_0}(\theta) =
\delta^q_{A_1}(\theta)$ for all~$\theta$.
\end{lem}

\begin{proof}
Let $(\S,\L)$ be a primitive $\Q$-concordance between two links $L_0$
and $L_1$.  Choose a primitive pair $(W_0,M_0)$ where
$\sigma_{k,d}(L_0)$ can be defined. Let
$(W,M)=(W_0,M_0)\cup_{(\Sigma_0,L_0)}(\S,\L)$.  Then $M$ is also
primitive.  We will compute $\sigma_{k,d}(L_1)$ from $(W,M)$.

Note that for any sufficiently large prime~$p$, the pair $(\S,\L)$ can
be considered as a concordance between links in $\Z_p$-homology
spheres and in particular $H_{q+1}(E_{\L},E_{L_0};\Z_p)=0$. For such a
prime $p$, let $(\tilde\S,\tilde\L)$ and $(\tilde\Sigma_0,\tilde L_0)$
be the $p$-fold branched covers of $(\S,\L)$ and $(\Sigma_0,L_0)$,
respectively.

We need the following consequence of Milnor's exact sequence for
infinite cyclic coverings~\cite{M1}. Suppose that $(\tilde X,\tilde
A)$ be a $p$-fold cyclic cover of a $CW$-complex $(X,A)$ induced by a
map $H_1(X)\to \Z_d$ that factors through the projection $\Z \to
\Z_d$. Then if $H_i(X,A;\Z_p)$ is zero, $H_i(\tilde X,\tilde A;\Z_p)$
is also zero.

In our case, $H_{q+1}(E_\L, E_{L_0};\Z_p)$ is zero and so
$H_{q+1}(E_{\tilde\L},E_{\tilde L_0};\Z_p)$ is also zero, for any
large prime $p$. From a Mayer-Vietoris sequence for
$(\tilde\S,\tilde\Sigma_0) = (E_{\tilde\L},E_{\tilde L_0})
\cup((\tilde\L,\tilde L_0)\times D^2)$, it follows that
$H_{q+1}(\tilde\S,\tilde\Sigma_0;\Z_p)$ is zero.  Hence $\tilde\S$
contributes nothing to the signature of the branched covering of
$(W,M)$, and $\sigma_{k,p}(W,M)=\sigma_{k,p}(W_0,M_0)$ by additivity
of signatures. Therefore $\sigma_{k,p}(L_0) = \sigma_{k,p}(L_1)$ and
by Lemma~\ref{lem:sign-branch-cover}, $\sigma^q_{A_0}(\theta) =
\sigma^q_{A_1}(\theta)$ if $\theta=2\pi k/p$ for some large prime~$p$
and some integer~$k$. Since the set of such $\theta$ is dense in the
real line, $\delta^q_{A_0}(\theta) = \delta^q_{A_1}(\theta)$ for
all~$\theta$.
\end{proof}

Now we are ready to show the invariance of the signature jump function
under $\Q$-concordance. The key to handle the non-primitive case is to
combine the geometric result for the primitive case with the
reparametrization formula in Section~\ref{sec:matrix-sign-poly}.

\begin{proof}[Proof of Theorem~\ref{thm:link-sign-inv}]
We may similarly proceed as in the proof of
Theorem~\ref{thm:conc-matrix}.  (Indeed, for the case of $q>1$, the
result is a consequence of Theorem~\ref{thm:conc-matrix}.) Let $L_0$
and $L_1$ be $\Q$-concordant links admitting the given type $\tau$.
Let $A_0$ and $A_1$ be Seifert matrices of Seifert surfaces of
$i_{c_0} L_0^\tau$ and $i_{c_1} L_1^\tau$, respectively. Let $\L^\tau$
be the union of $\tau_i$ parallel copies of the $i$-th component of a
$\Q$-concordance between $L_0$ and $L_1$. If $q=1$, we need to take
parallel copies with respect to the framing determined by $\tau$.  As
before, we can choose a positive integer $r$ such that $i_r \L^\tau$
is primitive.  Then $i_r i_{c_1} i_{c_0} L_0^\tau$ and $i_r i_{c_0}
i_{c_1} L_1^\tau$ are $\Q$-concordant via the primitive
$\Q$-concordance $i_r i_{c_0} i_{c_1} \L^\tau$.  Note that the links
$i_r i_{c_1} i_{c_0} L_0^\tau$ and $i_r i_{c_0} i_{c_1} L_1^\tau$ are
viewed as the boundaries of unions of parallel copies of Seifert
surfaces of the links $i_{c_0} L_0^\tau$ and $i_{c_1} L_1^\tau$. In
the case of $q=1$, the framing used to take $i_{c_k} L_k^\tau$ and the
framing induced by the Seifert surface of $i_{c_k} L_k^\tau$ must be
compatible in order for $i_r i_{c_1} i_{c_0} L_0^\tau$ and $i_r
i_{c_0} i_{c_1} L_1^\tau$ to be $\Q$-concordant. This is where our
framing condition in the definition of generalized Seifert surfaces
for 1-links is necessary.

The primitive links $i_r i_{c_1} i_{c_0} L_0^\tau$ and $i_r i_{c_0}
i_{c_1} L_1^\tau$ have the Seifert matrices $i^q_ri^q_{c_1}A_0$ and
$i^q_ri^q_{c_0}A_1$, respectively.  By Theorem~\ref{thm:sign-formula}
and Lemma~\ref{lem:prim-conc-sign}, we have
\begin{equation*}\begin{split}
\delta^\tau_{L_0}(\theta)
& = \sgn c_0\cdot \delta^q_{A_0}(\theta/c_0)\\
& = \sgn c_0\sgn c_1\cdot
\delta^q_{i^q_r i^q_{c_1} A_0}(\theta/r c_0 c_1) \\
& = \sgn c_0\sgn c_1\cdot
\delta^q_{i^q_r i^q_{c_0} A_1}(\theta/r c_0 c_1) \\
& = \sgn c_1\cdot \delta^q_{A_1}(\theta/c_1)
= \delta^\tau_{L_1}(\theta).
\end{split}\end{equation*}

This proves that $\delta^\tau_L(\theta)$ is a $\Q$-concordance
invariant independent of the choice of generalized Seifert surfaces.
\end{proof}

We remark that our framing condition for generalized Seifert surfaces
of 1-links seems to be the minimal requirement for signatures to be
well-defined $\Q$-concordance invariants.  For example, as an attempt
to remove the framing condition, one may fix corresponding framings on
$L_k$ for $k=0,1$ and consider $i_r L_k^\tau$ taken with respect to
the framings, where $r$ is a nonzero multiple of the complexity of
concordance. It can be shown that they admit Seifert surfaces, and are
$\Q$-concordant via a primitive concordance.  The signature jump
functions of Seifert matrices of $i_rL_k$ are equal by
Lemma~\ref{lem:prim-conc-sign}.  But the above argument can not be
applied to prove that the reparametrizations $\theta=\theta/r$ of the
signature jump functions are invariants of the original links $L_k$
and in fact they are not in general.

We also remark that it is meaningful to consider types for higher
dimensional links as well as 1-links. In fact the signature jump
function of type $\tau$ for a link $L$ is no more than the signature
jump function of type $(1\cdots 1)$ for the link $L^\tau$ obtained by
taking parallel copies, however, more information on $L$ is obtained
from types other than $(1\cdots 1)$ even for higher dimensional links
in the sphere. For example, suppose $K$ is a knot that is not torsion
in the (algebraic) knot concordance group (so that $\delta_K(\theta)$
is nontrivial~\cite{L2}) and $L$ be a split union of $K$ and
$-K$. Then $\delta^{(1,1)}_L$ is trivial, but if $a\ne b$,
$\delta^{(a,b)}_L(\theta)$ is nontrivial for some~$\theta$.

For $\tau=(1\cdots 1)$, we denote $\delta^\tau(L)$ and $c^\tau(L)$ by
$\delta(L)$ and $c(L)$, respectively.

Now we prove Theorem~\ref{thm:link-sign-prop}.

\begin{proof}[Proof of Theorem~\ref{thm:link-sign-prop}]

1. If $L$ is slice, $\delta^\tau_L$ is equal to the signature of the
unlink by Theorem~\ref{thm:link-sign-inv}. Since the null matrix is a
Seifert matrix for the unlink, $\delta^\tau_L(\theta)=0$.

2. Choose a positive multiple $r$ of $c^\tau (L)$, and let $A$ be a
Seifert matrix for $i_rL^\tau$. Then since $i^q_\alpha A$ is a Seifert
matrix of $i_r i_\alpha L^\tau$, $\delta_{i_\alpha L}^\tau(\theta) =
\delta^q_{i^q_\alpha A}(\theta/r) = \sgn n_\alpha\cdot
\delta^q_A(n_\alpha \theta/r) = \sgn
n_\alpha\cdot\delta_L^\tau(n_\alpha \theta)$ by
Theorem~\ref{thm:sign-formula}.

3. Let $A$ be a Seifert matrix of $i_{c^\tau(L)}L^\tau$. Since
$\delta^q_A(\theta)$ is of period $2\pi$, $\delta^\tau_L(\theta) =
\delta^q_A(\theta/c^\tau(L))$ is of period $2\pi\cdot c^\tau(L)$.

4.  For $i=1,2$, let $B_i$ be a $(2q+1)$-ball in the ambient space
$\Sigma_i$ of $L_i$ such that the connected sum $(\Sigma,L)$ is
obtained by gluing $(\Sigma_i-\inte(B_i),L_i-\inte(B_i))$ along their
boundaries. We view $\Sigma_i-\inte(B_i)$ as a subspace of $\Sigma$.
Let $c$ be a positive common multiple of $c^\tau(L_1)$ and
$c^\tau(L_2)$.  Then there is a generalized Seifert surface $F_i$ of
type $\tau$ with complexity $c$ for $L_i$. We may assume that each
component of $F_i \cap B_i$ is a $2q$-disk whose boundary contains
exactly one component of $\partial F_i \cap B_i $ so that the
intersection of $F_1-\inte(B_1)$ and $F_2-\inte(B_2)$ consists of
disjoint disks, and their union is a generalized Seifert surface $F$
of type $\tau$ with complexity $c$ for $L$ in $\Sigma$.

For $q>1$, let $A_i$ be a Seifert matrix of $F_i$. We have $H_q(F) =
H_q(F_1)\oplus H_q(F_2)$, and $A_1\oplus A_2$ is a Seifert matrix of
$F$.  Therefore $\delta^\tau_{L}(\theta) = \delta^q_{A_1\oplus
A_2}(\theta/c) = \delta^q_{A_1}(\theta/c) + \delta^q_{A_2}(\theta/c) =
\delta^\tau_{L_1}(\theta) + \delta^\tau_{L_2}(\theta)$.

For $q=1$, $L_1$ and $L_2$ are 1-knots by hypothesis.  We need
additional arguments since $H_1(F)$ is not decomposed as above. For a
manifold $E$, let denote the cokernel of $H_i(\partial E)\to H_i(E)$
by $\bar H_i(E)$. We assert that for any 1-knot and a generalized
Seifert surface $E$ of complexity $c$, the Seifert pairing on $H_1(E)$
induces a well-defined ``Seifert pairing'' on $\bar H_1(E)$ that
induces the signature jump function of the knot. For, if $x$ is a
boundary component of $E$, a parallel copy of $E$ is a 2-chain whose
boundary is homologous to $cx$ in the complement of $\inte(E)$ and so
the linking number of $x$ and any 1-cycle on $\inte(E)$ is zero. This
proves the assertion.

Assume that $F_i$ is connected for $i=1,2$.  Let $x_0,\ldots,x_{c-1}$
be boundary components of $F_1-\inte (B_1)$. $\{x_i\}$ are ordered so
that parallel cycles in $\Sigma-F$ are obtained by pushing $x_i$ and
$x_{i+1}$ slightly along the negative and positive normal directions
of $F$, respectively.  Let $y_1,\ldots,y_{c-1}$ be curves on
$\inte(F)$ such that the intersection number of $y_i$ and $x_j$ on $F$
is $\delta_{(i-1)j}-\delta_{ij}$, where $\delta_{ij}$ is the Kronecker
symbol (see Figure~\ref{fig:csum-generator} for a schematic picture).
Then $\bar H_1(F)$ is the direct sum of $\bar H_1(F_1)$, $\bar
H_1(F_2)$ and a free abelian group generated by $x_1,\ldots,x_{c-1}$
and $y_1,\ldots,y_{c-1}$.  We will show that the generators $x_i$ and
$y_i$ have no contribution to the signature.

Similarly to the proof of the previous assertion, it is shown that the
linking number of $x_i$ and any 1-cycle on the interior of
$F_j-\inte(B_j)$ is zero.  Furthermore, the value of the Seifert
pairing at $x_i$ and $y_j$ can be computed by counting intersections
of the 1-cycle obtained by pushing $y_i$ slightly and the 2-chain
obtained by attaching $(c-1)$ annuli bounded by $x_i-x_j$ ($j\ne i$)
to $F_1-\inte(B_1)$.  By this argument we can verify that the Seifert
matrix on $\bar H_1(F)$ with respect to generators of $\bar H_1(F_i)$
and $y_1,\ldots,y_{c-1}$, $x_1,\ldots,x_{c-1}$ is given by
$$
\arraycolsep=.3em
\def\fr#1#2{\frac{\mathstrut #1}{\mathstrut #2}}
A=\left[\begin{array}{cc|ccccc|ccccc}
A_1&\fr{}{}&*&*&*&\cdots&*& \\
\fr{}{}&A_2&*&*&*&\cdots&*& \\
\hline
*&*&*&*&*&\cdots&*& -\fr1c & -\fr1c & -\fr1c & \cdots & -\fr1c \\
*&*&*&*&*&\cdots&*& \fr{c-2}c & -\fr2c & -\fr2c & \cdots & -\fr2c \\
*&*&*&*&*&\cdots&*& \fr{c-3}c & \fr{c-3}c & -\fr3c & \cdots & -\fr3c \\
\vdots&\vdots&\vdots&\vdots&\vdots&\ddots&\vdots& \vdots & \vdots & \vdots & \ddots & \vdots \\
*&*&*&*&*&\cdots&*& \fr{1}c & \fr{1}c & \fr{1}c & \cdots & -\fr{c-1}c \\
\hline
& & \fr{c-1}c & \fr{c-2}c & \fr{c-3}c & \cdots & \fr1c & \\
& & -\fr1c & \fr{c-2}c & \fr{c-3}c & \cdots & \fr1c & \\
& & -\fr1c & -\fr2c & \fr{c-3}c & \cdots & \fr1c & \\
& &\vdots&\vdots&\vdots&\ddots&\vdots& \\
& & -\fr1c & -\fr2c & -\fr3c & \cdots & \fr1c & \\
\end{array}\right]
$$
where $A_i$ is a Seifert matrix defined on $\bar H_1(F_i)$ and empty
blocks represent zero matrices.  The bottom-center block and
middle-right block are nonsingular. In fact, subtracting the $i$-th
row from the $(i+1)$-th row for each $i=c-2,\ldots,1$, the
bottom-center block becomes an upper-triangular nonsingular matrix,
and similarly for the right-middle block. Hence the top-center block
and middle-left block can be eliminated by row and column
operations. This shows that $A$ is congruent to the block sum of
$A_1$, $A_2$ and a metabolic matrix.  Therefore $[A]=[A_1]+[A_2]$ in
the matrix cobordism group~$G^\Q_{+1}$, and we have $\delta_L(\theta)
= \delta^+_A(\theta/c) = \delta^+_{A_1}(\theta/c) +
\delta^+_{A_2}(\theta/c) = \delta_{L_1}(\theta) +
\delta_{L_2}(\theta)$ as desired.
\begin{figure}[tb]
\begin{center}
\leavevmode\includegraphics{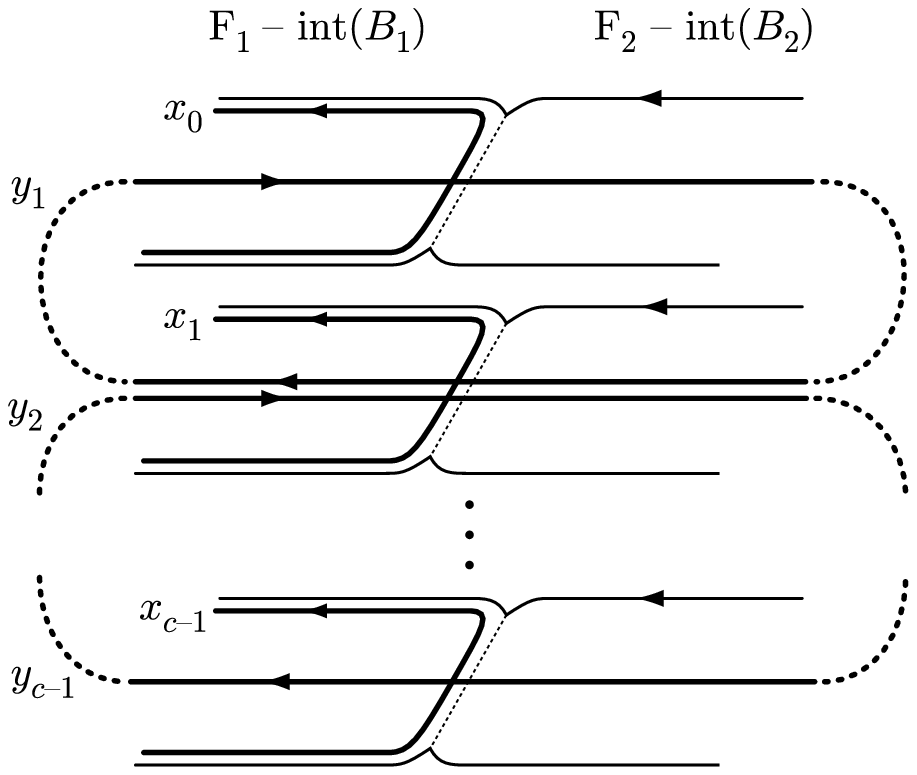}
\caption{}
\label{fig:csum-generator}
\end{center}
\end{figure}
\end{proof}

By Theorem~\ref{thm:link-sign-prop}, the signature jump function of a
link $L$ in an integral homology sphere always has the period $2\pi$
since $L^\tau$ is primitive.  However the period of the signature jump
function of a link with nontrivial complexity in a rational homology
sphere is not $2\pi$ in general. In the following example, for any
nonzero rational number~$r$ we construct a knot whose signature jump
function has the period~$2\pi r$.

\begin{exmp}
Let $K$ be a knot in $S^{2q+1}$ and $A$ be a Seifert matrix defined on
a Seifert surface $F$ of $K$. For $k\ne 0$, a knot $K'$ with Seifert
matrix $i^q_k A$ is obtained by fusing the link $i_k K=\partial(i_k
F)$ along $|k|-1$ bands whose interiors are disjoint to $i_k F$. (We
may take as $K'$ a $(k,1)$-cable knot with companion $K$ for $q=1$.)

Let $S$ be an embedded $(2q-1)$-sphere with knot type $K'$ and $J$ be
an unknotted $(2q-1)$-sphere that bounds a ball disjoint to $S$ in
$S^{2q+1}$. Let $C$ be a simple closed curve in $S^{2q+1}-S\cup J$
whose linking numbers with $S$ and $J$ are $n$ and $-1$, respectively,
as shown in Figure~\ref{fig:period-knot}. For $n\ne 0$, the result of
surgery along $S$ and $C$ (with respect to zero-framing if $q=1$) is a
rational homology sphere~$\Sigma$, and $J$ can be considered as a knot
in~$\Sigma$.

\begin{figure}[hbt]
\begin{center}
\leavevmode\includegraphics{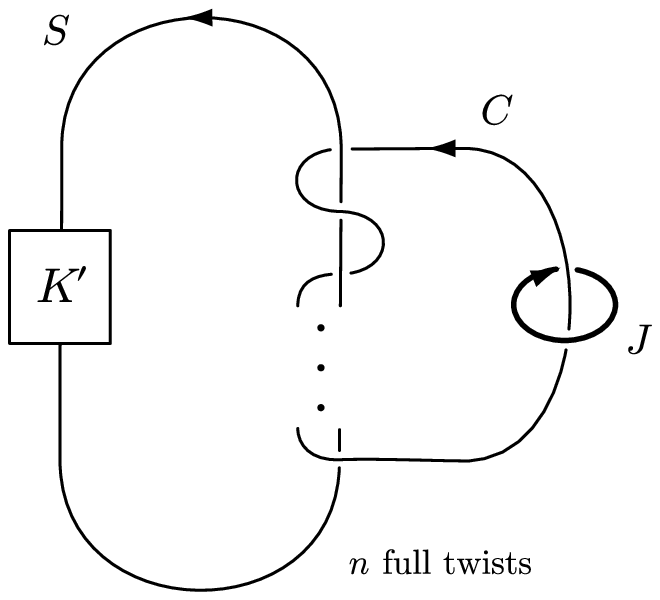}
\caption{}
\label{fig:period-knot}
\end{center}
\end{figure}

$H_1(\Sigma-J)$ is an infinite cyclic group generated by a meridian
$u$ of~$S$. The meridian $\mu$ of $J$ is given by $\mu=nu$ in
$H_1(\Sigma-J)$.  Hence $c(J)=|n|$ by
Theorem~\ref{thm:compute-complexity}. In fact, there is a submanifold
in $S^{2q+1}$ bounded by the union of a parallel of $S$ and $i_n J$ as
in Figure~\ref{fig:period-surface}, and a generalized Seifert surface
$E$ with complexity $n$ for $J$ is obtained by gluing a $2q$-disk in
$\Sigma$ to the submanifold along the parallel of~$S$.  For $q\ne 1$,
$i^q_k A$ is a Seifert matrix of~$E$. For $q=1$, the block sum of
$i^q_k A$ and a zero matrix of dimension $|n|-1$ is a Seifert matrix
of~$E$.  In any case, we have $\delta_J(\theta) = \sgn n\cdot
\delta^q_{i^q_k A}(\theta/n) = \sgn(k/n)\cdot \delta_K(\theta\cdot
k/n)$.
\begin{figure}[hbt]
\begin{center}
\leavevmode\includegraphics{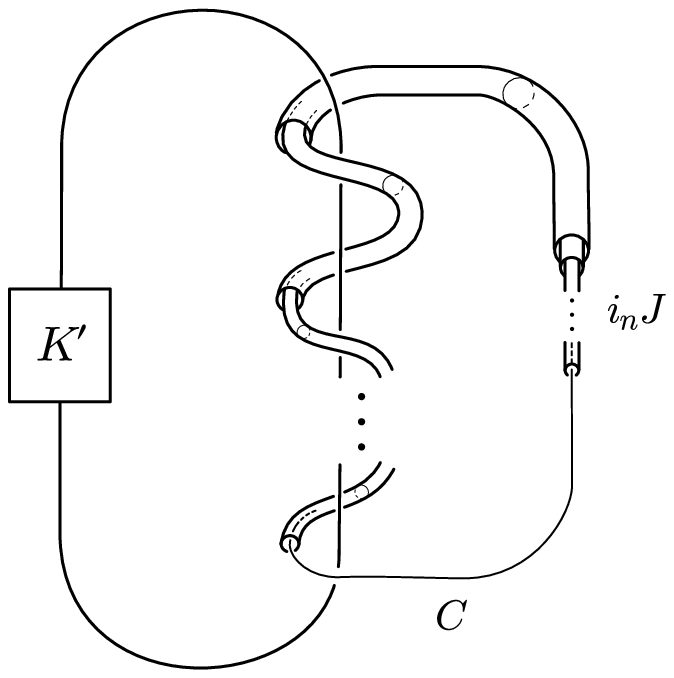}
\caption{}
\label{fig:period-surface}
\end{center}
\end{figure}
\end{exmp}

As a special case of this example, there are knots in rational
homology spheres whose signature jump functions do not have the period
$2\pi$. Such knots are not $\Q$-concordant to knots in (homology)
spheres, by Theorem~\ref{thm:link-sign-prop}. This illustrates that
the $\Q$-concordance theory of links in rational homology spheres is
not reduced to that of links in spheres.

Sophisticating the above arguments, we show
Theorem~\ref{thm:independent-knots}.

\begin{proof}[Proof of Theorem~\ref{thm:independent-knots}]
Choose infinitely many distinct primes $p_i$ greater than~7.  For odd $q$,
let
$$
A=\begin{bmatrix}1 & 1 \\ 0 & 1\end{bmatrix}
$$
and for even $q$, let
$$
A=\begin{bmatrix}
1 & 1 & 0 & 0 \\ 0 & 0 & 1 & 0 \\ 0 & -1 & 0 & 1 \\ 0 & 0 & 0 & 1
\end{bmatrix}.
$$
It is easily checked that $\delta^q_A(\theta)$ is nonzero if and only
if $\theta=n\pi(k\pm 1/6)$ for some integer $k$, where $n=1$ for
even~$q$ or $n=2$ for odd~$q$.  Since $A$ is a Seifert matrix of a
knot in $S^{2q+1}$ by~\cite{Ker,L1}, a $(2q-1)$-knot $K_i$ in a
rational homology sphere such that $\delta_{K_i}(\theta) =
\delta^q_A(\theta\cdot p_i/7)$ is constructed by the above example.
(For $q=2$, we need to use the block sum of 8 copies of $A$ to satisfy
the extra condition~\cite{L1} forced by a result of Rochlin. Since the
signature jump function multiplies by 8, this is irrelevant to our
purpose.) In order to show that $K_i$'s are independent, we consider a
connected sum $K$ of $a_i K_i$ over all~$i$, where all but finitely
many $a_i$ are zero and $a_j\ne 0$ for some $j$.  We will show that
$K$ is not $\Q$-concordant to a knot in a integral homology sphere.
This will complete the proof.

By the additivity of signature jump function, $\delta_K(\theta)=\sum_i
a_i \delta^q_A(\theta\cdot p_i/7)$.  We claim that
$\delta_K(7n\pi/6p_j) \ne 0$ but $\delta_K(7n\pi/6p_j+2\pi) = 0$. If a
summand $a_i\delta^q_A(n\pi p_i/6p_j)$ of $\delta_K(7n\pi/6p_j)$ is
nonzero, then $n\pi p_i/6p_j = n\pi(k\pm 1/6)$ for some $k$, or
equivalently $p_i = p_j(6k\pm 1)$. This implies $p_j \mid p_i$ and
$i=j$. This shows the first claim. On the other hand, if a summand
$a_i\delta^q_A(n\pi p_i/6p_j+2\pi p_i/7)$ of
$\delta_K(7n\pi/6p_j+2\pi)$ is nonzero, then we have
$p_i(7n+12p_j)=7np_j(6k\pm 1)$ for some $k$, and $7\mid p_j$ or $7\mid
p_j$. This is a contradiction, and the second claim is proved.

From the claims, $\delta_K(\theta)$ is not of period $2\pi$. By
Theorem~\ref{thm:link-sign-prop}, $K$ is not $\Q$-concordant to knots
of period $2\pi$, in particular, to knots in integral homology
spheres.
\end{proof}

Another natural question is whether the homomorphism
$C^{\vphantom{\Q}}_n \to C^\Q_n$ is injective. More generally, one may
ask whether the natural map from the set of concordance classes of
links in spheres into that of links in rational homology spheres is
injective. We do not address these questions in this paper.

We finish this paper with two applications illustrating that our
results are also useful for links in spheres.  One is a simple proof
of the following result of Litherland~\cite{Li}: Let $K$ be a
satellite knot with companion $J$, that is, $K = h(J)$ where $J$ is a
1-knot in a standard solid torus $V$ in $S^3$ and $h$ is a
homeomorphism from $V$ onto a tubular neighborhood of a knot $K'$
sending a 0-linking longitude on $\partial V$ to a 0-linking longitude
of $K'$. Let $r$ be the winding number of $J$ in $V$.  Then
$\delta_{K}(\theta) = \delta_J(\theta)+\delta_{K'}(r\theta)$. For, if
$A$ and $B$ are Seifert matrices of $J$ and $K'$, then the block sum
$A\oplus i^+_\alpha B$ is a Seifert matrix of $K$ for some tuple
$\alpha$ such that $n_\alpha=r$, and Litherland's formula is proved by
Theorem~\ref{thm:sign-formula}.  We remark that a formula for
Alexander polynomials of satellite knots~\cite{Se} can also be proved
using an analogous result of Theorem~\ref{thm:sign-formula} mentioned
in Section~\ref{sec:matrix-sign-poly}.

Another is a new proof of (a weaker version of) a result of
Kawauchi~\cite{Ka}: If a $(2q-1)$-knot $K$ is not (algebraically if
$q=1$) torsion, then $i_r K$ is a boundary link that is never
concordant to any split union for any $r>0$. For, if $i_rK$ is
concordant to a split link, $i_2K$ is concordant to a split union
$K\cup K$, and then by Theorem~\ref{thm:link-sign-inv}
and~\ref{thm:link-sign-prop}, $\delta_K(2\theta) =
\delta_{i_2K}(\theta) = \delta_{K\cup K}(\theta) =
2\delta_K(\theta)$. This implies that $\delta_K(\theta)=0$ for all
$\theta$ and $K$ is (algebraically if $q=1$) torsion by a result of
Levine~\cite{L2}.


\begin{thebibliography}{00}


%\bibitem{AK} S.~Akbulut and R.~Kirby,
%{\em Branched coverings of surfaces in 4-manifolds},
%Math.\ Annalen (1980), 111--131.

%\bibitem{CG2} A.~Casson and C.~Gordon,
%{\em On slice knots in dimension three},
%Proc.\ Symp.\ in Pure Math.\ XXX (1978), part two, 39--53.

\bibitem{CK1} J.~Cha and K.~Ko,
\emph{On equivariant slice knots},
Proc.\ A.~M.~S., 127 (1999), 2175--2182.

%\bibitem{CK2} J.~Cha and K.~Ko,
%\emph{Signature invariants of links from irregular covers and non-abelian
%covers},
%Math.\ Proc.\ Camb.\ Phil.\ Soc.\ 127 (1999), 67--81.

\bibitem{CK4} J.~Cha and K.~Ko,
\emph{Signatures of covering links},
arXiv:math.GT/0108206.

%\bibitem{CL} T.~Cochran and J.~Levine,
%{\em Homology boundary links and the Andrew-Curtis conjecture},
%Topology 29 (1990), 231--239.

\bibitem{CO1} T.~Cochran and K.~Orr,
{\em Not all links are concordant to boundary links},
Bull.\ A.~M.~S. 23 (1990), 99--106.

\bibitem{CO2} T.~Cochran and K.~Orr,
{\em Not all links are concordant to boundary links},
Ann.\ Math.\ 138 (1993), 519--554.

%\bibitem{CO3} T.~Cochran and K.~Orr,
%{\em Homology boundary links and Blanchfield forms: concordance
%classification and new tangle-theoretic construction},
%Topology 33 (1994), 397--427.

%\bibitem{CO4} T.~Cochran and K.~Orr,
%\emph{Homology cobordism and generalizations of Milnor's invariants},
%J. of Knot Theory and Its Ramifications 8 (1999), 429--436.

%\bibitem{CS} S.~Cappell and J.~Shaneson,
%\emph{Link Cobordism},
%Comm.\ Math.\ Helv.\ 55 (1980), 20--49.

\bibitem{Er} D.~Erle,
\emph{Die quadratische Form eines Knottens und ein Satz \"uber
Knotten-mannigfaltigkeiten},
J.~Reine Angw.\ Math.\ 236 (1969), 174--218.

\bibitem{G} P.~Gilmer,
\emph{Link cobordism in rational homology 3-spheres},
J.\ Knot Theory and Its Ramifications 2 (1993), no 3, 285--320.

\bibitem{GL} P.~Gilmer and C.~Livingston,
{\em The Casson-Gordon invariant and link concordance},
Topology 31 (1992), 475--492.

\bibitem{Hi} J.~Hillman,
\emph{Alexander ideals of links},
Springer-Verlag Lecture Notes in Math.\ 895, Springer, Berlin, 1981.

\bibitem{K} L.~Kauffman,
{\em On Knots}, Princeton University Press, 1987.

\bibitem{KT} L.~Kauffman and L.~Taylor,
{\em Signatures of links},
Trans.\ A.~M.~S. 216 (1976), 351--365.

\bibitem{Ka} A.~Kawauchi,
\emph{On links not cobordant to split links},
Topology 19 (1980), 321--334.

\bibitem{Ke} C.~Kearton,
\emph{Blanchfield duality and simple knots},
Trans.\ A.~M.~S. 202 (1975), 141--160.

\bibitem{Ker} M.~Kervaire,
\emph{Les noeuds de dimensions sup\'erieres},
Bull.\ Soc.\ Math.\ France 93 (1965), 225--271.

\bibitem{Ko} K.~Ko,
{\em Seifert matrices and boundary link cobordisms},
Trans.\ A.~M.~S. 299 (1987), 657--681.

\bibitem{L0} J.~Levine,
{\em Unknotting spheres in codimension two},
Topology 4 (1965), 9--16.

\bibitem{L1} J.~Levine,
{\em Knot cobordism group in codimension two},
Comm.\ Math.\ Helv.\ 44 (1969), 229--244.

\bibitem{L2} J.~Levine,
{\em Invariants of knot cobordism},
Invent.\ Math.\ 8 (1969), 98--110.

\bibitem{L4} J.~Levine,
\emph{Link invariants via the eta invariant},
Comment.\ Math.\ Helv.\ 69 (1994), 82--119.

\bibitem{Li} R.~Litherland,
\emph{Signatures of iterated torus knots},
Topology of low-dimensional manifolds
(Proc. Second Sussex Conf., Chelwood Gate, 1977), pp.\ 71--84,
Lecture Notes in Math., 722, Springer, Berlin, 1979.

\bibitem{Ma} T.~Matumoto,
{\em On the signature invariants of a non-singular complex sesqui-linear form},
J. Math.\ Soc.\ Japan 29 (1977), 67--71.

\bibitem{M1} J. Milnor,
{\em Infinite cyclic coverings},
Conference on the Topology of Manifolds,
Prindle, Weber \& Schmidt, Boston, Mass., 1968, pp. 115--133.

%\bibitem{Mi} W.~Mio,
%\emph{On boundary-link cobordism},
%Math.\ Proc.\ Camb.\ Phil.\ Soc.\ 101 (1987), 542--560.

\bibitem{Mu} K.~Murasugi,
{\em On a numerical invariant of link types},
Trans.\ A.~M.~S.\ 117 (1965), 387--422.

\bibitem{Se} H.~Seifert,
\emph{On the homology invariants on knots},
Quart.\ J. Math.\ Oxford 2 (1950), 23--32.

\bibitem{Tr} A.~G.~Tristram,
\emph{Some cobordism invariants for links},
Proc.\ Camb.\ Phil.\ Soc.\ 66 (1969), 252--264.

\bibitem{Tro} H.~F.~Trotter,
\emph{Homology of group systems with applications to knot theory},
Ann.\ Math.\ 59 (1962), 464--498.

\end{thebibliography}
\end{document}